\documentclass[preprint]{elsarticle}  

\usepackage[linesnumbered,vlined, procnumbered,titlenumbered,ruled]{algorithm2e}
 
\usepackage{lipsum}
\makeatletter
\def\ps@pprintTitle{%
 \let\@oddhead\@empty
 \let\@evenhead\@empty
 \def\@oddfoot{}%
 \let\@evenfoot\@oddfoot}
\makeatother 
 
 \usepackage{caption}
\usepackage{setspace}
\usepackage{amssymb}
 \usepackage{amsthm}
\usepackage{amsmath}
\usepackage{tikz}
\usepackage{hhline}
\usetikzlibrary{decorations.pathreplacing}
\usetikzlibrary{decorations.pathmorphing}

\theoremstyle{plain}
\newtheorem{theorem}{Theorem}[section]

\newtheorem{lemma}[theorem]{Lemma}

\journal{--}

\begin{document}

\begin{frontmatter}

\title{The Enumeration of Cyclic MNOLS}

\author{
Fatih Demirkale
}
\ead{
fatihd@yildiz.edu.tr
}
\address{Department of Mathematics \\ Y{\i}ld{\i}z Technical University\\
Esenler, \.{I}stanbul 34220, Turkey\\
}

\author{
Diane M. Donovan
}
\ead{dmd@maths.uq.edu.au}
\address{
School of Mathematics and Physics\\ 
Centre for Discrete Mathematics and Computing\\ 
University of Queensland\\
St Lucia 4072, Australia \\
}

\author{
Janne I. Kokkala\fnref{label4}
}
\fntext[label4]{This work was supported
in part by the Academy of Finland, Project \#289002.}
\ead{janne.kokkala@iki.fi}
\address{
Department of Communications and Networking\\
School of Electrical Engineering, P.O. Box 13000, 00076\\
Aalto University\\
Aalto, Finland\\
}

\author{Trent G. Marbach\fnref{label2}} %
\fntext[label2]{This research was supported in part by the Monash eResearch Centre and eSolutions-Research Support Services through the use of the MonARCH HPC Cluster.}
\ead{trent.marbach@uqconenct.edu.au}
\address{School of Mathematics and Physics\\
University of Queensland\\
St Lucia 4072, Australia}
\address{School of Mathematical Sciences\\
Monash University\\
Clayton  3800, Australia}

\doublespacing 
\begin{abstract}
In this paper we study collections of mutually nearly orthogonal Latin squares ($\text{MNOLS}$), which come from a modification of the orthogonal condition for  mutually orthogonal Latin squares. 
 In particular, we find the maximum $\mu$ such that there exists a set of $\mu$ cyclic $\text{MNOLS}$ of order $n$ for $n \leq 18$, as well as providing a full enumeration of sets and lists of $\mu$ cyclic $\text{MNOLS}$ of order $n$ under a variety of equivalences with $n \leq 18$. 
This resolves in the negative a conjecture that proposed the maximum $\mu$ for which a set of $\mu$ cyclic $\text{MNOLS}$ of order $n$ exists is $\lceil n/4\rceil +1$.

\end{abstract}


\end{frontmatter}

\doublespacing 
\section{Introduction}

The study of mutually orthogonal Latin squares ($\text{MOLS}$) is a subject that has attracted much attention. 
Such interest has been stimulated by the relevance of the field, with applications in error correcting codes, cryptographic systems, affine planes, compiler testing, and statistics (see \cite{DisMathUsingLS}). 
Although, as it is well known, there exists a set of $n-1$ $\text{MOLS}$ of order $n$ when $n$ is a prime or a prime power, the largest number of $\text{MOLS}$ of order $n$ known to exist  when $n$ is even is generally much smaller and such sets of $\text{MOLS}$ are hard or impossible to find; there does not exist $2$ $\text{MOLS}$ of order $6$, and it is unknown whether three $\text{MOLS}$ of order $10$ exists or not. 

Formally, a \emph{Latin square} of order $n$ is an $n \times n$ array in which the  $n$ distinct symbols $\{0, \ldots, n-1\}$ are arranged so that each symbol occurs once in each row and once in each column. 
We index the rows and columns by $\{0, \ldots, n-1\}$. 
For a Latin square $L$ of order $n$ we may write $(r,c,e)\in L$ to mean that $L$ has a cell in row $r$ and column $c$ that contains symbol $e$.   
We assume the 3 components of such a triple are taken mod $n$, so that $(r,c,e+n) = (r,c+n,e)=(r+n,c,e)$. 
The notation $(r,c,e)\in L$ is known as orthogonal array notation. 
We also write $L(r,c)=e$ when $(r,c,e) \in L$. 
A pair of Latin squares $L_1, L_2$ of order $n$ are called \emph{orthogonal} if the superimposition of $L_1$ and $L_2$ contains each ordered pair of symbols exactly once. 
A set of $\mu$ Latin squares are \emph{mutually orthogonal} if they are pairwise orthogonal, and we refer to such a set as a \emph{set of  $\mathrm{MOLS}$}.

Based upon the significance and usefulness that is exhibited in the study of $\text{MOLS}$, Raghavarao, Shrikhande, and Shrikhande \cite{MNOLSinvent} introduced a modification to the definition of orthogonality to overcome restrictions for the even order case. 
A pair of Latin squares $L_1, L_2$ of even order $n$ are called \emph{nearly orthogonal} if the superimposition of $L_1$ and $L_2$ contains each ordered pair of symbols $(\ell,\ell')$ exactly once, except in the case $\ell = \ell'$, where no such pair occurs, and in the case $\ell \equiv \ell' +n/2 \pmod{n}$, where such pairs occur twice. 
We consider collections of $\mu$ Latin squares of order $n$ that are pairwise nearly orthogonal, which are denoted as collections of $\mu$ $\text{MNOLS}$ of order $n$. 
Traditionally these collections are unordered sets, although we will also consider ordered lists.

An \emph{orderly algorithm} is a way of generating all examples of some combinatorial object, such that all equivalence classes appear in the generation, but during the generation no two objects constructed are equivalent. 
This technique is typically attributed to \cite{Faradzev} and \cite{EveryOne}. 
A similar technique, called canonical augmentation \cite{IsoFreeExhaGen},  has been used to generate Latin rectangles by augmenting a row at a time (see also \cite{LS8}\cite{MostLSsubsquares}).
This is not the only method of enumerating Latin rectangles, and a variety of enumerative techniques have been applied to solve it (see \cite{numberLS} and the citations contained within).
Recently, this work has lead to the enumeration of  $\text{MOLS}$ for order less than or equal to $9$ \cite{EnumMOLS}. 
In a similar vein, we will perform an orderly algorithms that generate collections of $\mu$ cyclic $\text{MNOLS}$ of order $n$ under certain equivalences.  
See \cite{ClassAlgorithms} for a general reference on this kind of enumeration problem.

The pioneering work \cite{MNOLSinvent} on sets of $\mu$ $\text{MNOLS}$ of order $n$ investigated an upper bound on $\mu$ when $n$ is fixed, and they showed that if there exists a set of $\mu$ $\text{MNOLS}$ of order $n$, then $\mu \leq n/2 +1$ for $n \equiv 2 \pmod{4}$ and $\mu \leq n/2 $ for $n \equiv 0 \pmod{4}$. 
In the case that a set of $\mu$ $\text{MNOLS}$ of order $n$ obtains this bound, it is called a complete set of $\mu$ $\text{MNOLS}$ of order $n$. 
The authors proceeded to explore the existence of sets of $\mu$ $\text{MNOLS}$ of order $n$ by investigating sets of $\mu$ cyclic $\text{MNOLS}$ of order $n$; that is each Latin square $L$ has $L(r,c+1) \equiv L(r,c)+1 \pmod{n}$ for all $r,c \in [0,n-1]$, recalling that the entries are taken mod $n$. 
The sets of $\mu$ $\text{MNOLS}$ of order $n$ that were found included single examples of sets of two cyclic $\text{MNOLS}$ of order $4$, sets of three cyclic $\text{MNOLS}$ of order $6$, and sets of three cyclic $\text{MNOLS}$ of order $8$, demonstrating that the bound is tight for $n = 4$. 
It was later shown \cite{no4MNOLS6} that there does not exist a set of four $\text{MNOLS}$ of order $6$, and so the bound is not tight for $n=6$. 

Further results \cite{NOLS} showed sets of three $\text{MNOLS}$ of order $n$ exist for even $n \geq 358$. 
The authors also introduced a concept of equivalence between sets of $\mu$ cyclic  $\text{MNOLS}$ of order $n$ called isotopic equivalence (details  in Section \ref{sec2}). 
They found a number of isotopically non-equivalent sets of $\mu$ cyclic  $\text{MNOLS}$ of order $n$ for $n \leq 12$. 
The number of these sets of $\mu$ cyclic  $\text{MNOLS}$ of order $n$ is given in Table \ref{tab:NOLS}.

\begin{table}[!h]
\begin{center}
  \begin{tabular}{ |c | c | c|c | c | c|c | c | c|c | c | c|c | c | c|c | c | c|c|  }
    \hline
$n$& $6$& $8$& $10$& $12$ \\ \hline
$\mu=3$&1&1&$\geq1$&$>1$\\ \hline
$\mu=4$&0&0&1&$>1$ \\ \hline
$\mu=5$&0&0&0&$0$ \\ \hline

\end{tabular}
\end{center}
\caption{The number of sets of $\mu$ $\text{MNOLS}$ of order $n$ under isotopic equivalence.}\label{tab:NOLS}
\end{table}

\begin{table}[!h]
\begin{center}
  \begin{tabular}{ |c | c | c|c | c | c|c | c | c|c | c | c|c | c | c|c | c | c|c|  }
    \hline
$0$& $1$& $2$& $3$ \\ \hline
$1$& $2$& $3$& $0$ \\ \hline
$2$& $3$& $0$& $1$ \\ \hline
$3$& $0$& $1$& $2$ \\ \hline
\end{tabular}
\hspace{10pt}
  \begin{tabular}{ |c | c | c|c | c | c|c | c | c|c | c | c|c | c | c|c | c | c|c|  }
    \hline
$1$& $2$& $3$& $0$ \\ \hline
$3$& $0$& $1$& $2$ \\ \hline
$0$& $1$& $2$& $3$ \\ \hline
$2$& $3$& $0$& $1$ \\ \hline
\end{tabular}
\end{center}
\caption*{A pair of Latin squares of order $4$ that are nearly orthogonal.}
\end{table}

The literature surrounding the search for the existence of sets of three cyclic $\text{MNOLS}$ has been documented in \cite{DCA}. This document also gave further constructions that verified the existence of a set of three $\text{MNOLS}$ of order $n$ for all even $n \geq 6$, except perhaps when  $n= 146$. 

In the current paper, we find the maximum $\mu$ such that there exists a set of $\mu$ cyclic  $\text{MNOLS}$ of order $n$ for $n \leq 18$, as well as providing a full enumeration of sets and lists of $\mu$ cyclic  $\text{MNOLS}$ of order $n$ under a variety of equivalences with $n \leq 18$. 
This will resolve in the negative a conjecture of \cite{NOLS} that proposed the maximum $\mu$ for which a set of $\mu$ cyclic  $\text{MNOLS}$ of order $n$ exists is $\lceil n/4\rceil +1$ (the maximum $\mu$ appears erroneously as $\lceil n/8\rceil +1$ in the original conjecture \cite{PrivateVanRees}, and the maximum value we have written was the intended conjecture).

\section{Further definitions} \label{sec2}

A $(\mu,n)$-difference set is a set of $n$ $\mu$-tuples $\{(a_k^1, \ldots, a_k^{\mu}) \mid 1 \leq k \leq n\}$ over the alphabet $\{0, \ldots, n-1\}$ in which $a_k^i \neq a_l^i$ for all $1\leq k,l \leq n$ and $1 \leq i \leq \mu$, and  the multiset of ordered differences modulo $n$ between elements in two positions $i,j$, i.e.\ $\{a_k^i - a_k^j \pmod{n} \mid 1 \leq k \leq n \}$,
 does not contain $0$, contains $n/2$ twice, and contains every other difference once \cite{NOLS}. 
We can define $\mu$ Latin squares, $A_i$, that form a collection of $\mu$ cyclic $\text{MNOLS}$ of order $n$ from a $(\mu,n)$-difference set by defining the cells of the first columns as $A_i(r,0)=a_r^i$ and each cell in subsequent columns by adding $1\pmod{n}$ to the symbol of the corresponding cell in the previous column. It is also clear that this process is reversible, so a list of $\mu$ cyclic $\text{MNOLS}$ of order $n$ can be used to construct a $(\mu,n)$-difference set. 

We will enumerate both ordered lists and unordered sets of $\mu$ cyclic $\text{MNOLS}$ of order $n$. 
A \emph{set of $\mu$ $\mathrm{MNOLS}$} of order $n$ is a set $\{L_1, \ldots, L_{\mu}\}$  such that $L_i, L_j$ are nearly orthogonal for $1 \leq i, j \leq \mu$, $i \neq j$.  
A \emph{list of $\mu$ $\mathrm{MNOLS}$} of order $n$ is an ordered list $(L_1, \ldots, L_{\mu})$ such that $L_i, L_j$ are nearly orthogonal for $1 \leq i, j \leq \mu$, $i \neq j$.  
This distinction will be important when we enumerate collections of $\mu$ $\text{MNOLS}$ of order $n$. 
We will write \emph{collection} when a statement holds for either a list or set. 
A list of $\mu$ $\text{MNOLS}$ of order $n$ $(L_1, \ldots, L_{\mu})$ is \emph{reduced} if $L_1$ has its first row and column in natural order.

\section{Cyclic MNOLS}

We saw previously that a collection of $\mu$ cyclic $\text{MNOLS}$ of order $n$ can be developed from a $(\mu,n)$-difference set. 
Each $(\mu,n)$-difference set will correspond to $n!$ distinct lists of $\mu$ cyclic $\text{MNOLS}$ of order $n$, depending on which cells of the first column of the list of $\mu$ cyclic $\text{MNOLS}$ are chosen to be filled with which $\mu$-tuple from the $(\mu,n)$-difference set. 
This partitions the lists of $\mu$ cyclic $\text{MNOLS}$ of order $n$ evenly into classes of size $n!$, and so from now we will enumerate the number of $(\mu,n)$-difference sets under a number of equivalences, and consequentially we will have enumerated the number of collections of $\mu$ cyclic $\text{MNOLS}$ of order $n$.

Given a $(\mu,n)$-difference set, we may: 1) rearrange the order of all $\mu$-tuples simultaneously using the same rearrangement on all $\mu$-tuples; 
2) multiply all symbols in all $\mu$-tuples of the difference set by $a$ where $\gcd(a, n) = 1$; 
3) add a constant to all symbols in all $\mu$-tuples of the difference set; or any combination of 1), 2), and 3). These operations form a group $G$ that acts on the set of all $(\mu,n)$-difference sets. Let $G_1$ be the subgroup of operations of the form 1), $G_2$ be the subgroup of operations of the form 2), and $G_3$ be the subgroup of operations of the form 3).
We can think of the elements of $G$ as group actions that take a $(\mu,n)$-difference set and a group element and produce another $(\mu,n)$-difference set.

Two $(\mu,n)$-difference sets are \emph{set-isotopic} if one can be obtained from the other using the group actions in $G$, \emph{list-isotopic} if one can be obtained from the other using group actions in both $G_2$ and $G_3$, \emph{set-reduced-equivalent} if one can be obtained from the other using group actions in $G_1$, and \emph{list-reduced-equivalent} if they are identical as sets.

For a given $(\mu,n)$-difference set, $\mathcal{D}$, the set of all $(\mu,n)$-difference sets set-isotopic to $\mathcal{D}$ is called its \emph{set-isotopy class}. We similarly define the \emph{list-isotopy class}, \emph{set-reduced class}, and \emph{list-reduced class}. 
We let $\mathcal{C}^{\mu}_n$ denote a set of $(\mu,n)$-difference sets that contains precisely one representative from each set-isotopy class.

We define $\mathrm{Is}_s(\mathcal{L})$ to be the set of group actions of $G$ that fixes $\mathcal{L}$. 
Similarly we define $\mathrm{Is}_l(\mathcal{L})$, $\mathrm{Red}_s(\mathcal{L})$, and $\mathrm{Red}_l(\mathcal{L})$ respectively.

Two sets of $\mu$ $\text{MNOLS}$ of order $n$ are \emph{isotopic} if permuting the rows, columns, and symbols consistently among all Latin squares in one set yields the other set. 

\begin{lemma}
Two sets of $\mu$ cyclic $\text{MNOLS}$ of order $n$ are isotopic if and only if their two corresponding $(\mu,n)$-difference sets are also set-isotopic. 
\end{lemma}
\begin{proof}
The reverse direction is trivial, as we can change the rows and symbols between the $\text{MNOLS}$ in the same way that the rows and symbols were changed between the  $(\mu,n)$-difference sets. Some columns may need to be swapped, but this is easily seen.  

For the forward direction, let the two Latin squares be $L=[\ell_{ij}]$ and $L' = [\ell'_{ij}]$. Consider that the difference $\ell'_{ij} - \ell_{ij}$ is the same for each $j$ when we fix $i$. 
 This means a permutation of the symbols will send the pairs of symbols with difference $1$ to pairs of symbols of some constant difference $x$.  The only way this can happen is if the permutation is of the form $i \mapsto xi+j$ with $\gcd(x,n)=1$. A permutation of the rows will have no effect on the corresponding  $(\mu,n)$-difference sets, and any permutation of the columns is fixed by the permutation of the symbols (except for trivial cyclic shifts, which can be counteracted by changing the permutation of the symbols).  
\end{proof}

We call a pair of collections of $\mu$ cyclic $\text{MNOLS}$ of order $n$ set-isotopic (resp.\ list-isotopic, set-reduced-equivalent, list-reduced-equivalent) if their corresponding  $(\mu,n)$-difference sets are also set-isotopic (resp.\ list-isotopic, set-reduced-equivalent, list-reduced-equivalent). Note that this definition applies only to cyclic $\text{MNOLS}$ and not general $\text{MNOLS}$.

The computation in this paper will find the number of $(\mu,n)$-difference sets distinct up to list-isotopy, set-isotopy, list-reduction, and set-reduction for $n \leq 18$ and for each $2 \leq \mu \leq 5$.

\begin{lemma} \label{numberClassesPerClass}
Given $\mathcal{L} \in \mathcal{C}_n^{\mu}$:
\begin{enumerate}
\item the number of list-isotopy classes within the set-isotopy class of $\mathcal{L}$ is: \\
$ \mu! \cdot {|\mathrm{Is}_l(\mathcal{L})|}/{|\mathrm{Is}_s(\mathcal{L})|};$
\item the number of set-reduced classes within the set-isotopy class of $\mathcal{L}$ is: \\
$\phi(n) \cdot n \cdot {|\mathrm{Red}_s(\mathcal{L})|}/{|\mathrm{Is}_s(\mathcal{L})|};
\text{ and}$
\item the number of list-reduced classes within the set-isotopy class of $\mathcal{L}$ is: \\
$\phi(n) \cdot n \cdot \mu! \cdot {|\mathrm{Red}_l(\mathcal{L})|}/{|\mathrm{Is}_s(\mathcal{L})|}={\phi(n) \cdot n \cdot \mu!}/{|\mathrm{Is}_s(\mathcal{L})|}.$

\end{enumerate}
\end{lemma} 
\begin{proof}
By the orbit-stabilizer theorem. 
\end{proof}

\section{Canonical forms}
Given a partition of $\mathcal{C}^{\mu}_n$ as $\mathcal{C}^{\mu}_n = \cup_{i=1}^{\alpha} C_i$ with $C_i \cap C_j =\emptyset$ for $1 \leq i < j \leq \alpha$, a \emph{canonical form} is a function $f:\mathcal{C}^{\mu}_n \rightarrow \mathcal{C}^{\mu}_n$ such that for all $\mathcal{L},\mathcal{M} \in C_i$, $f(\mathcal{L}) = f(\mathcal{M})$ and $f(\mathcal{L})\in C_i$. 
We will say the lists within $\text{Im}(f)$ are \emph{canonical}. 
This allows us to represent each set-isotopy class of $(\mu,n)$-difference sets  by a single $(\mu,n)$-difference sets. 
For Latin squares, there are procedures to find a canonical Latin square amongst an isotopy class, which usually make use of a program to find a canonical labelling of a graph, with implementations such as nauty \cite{nauty}. 
While this is usually a computationally intensive task (having no known polynomial time algorithm), it will not form the bottleneck of the search, and so there is no need for such a sophisticated method. 
As such, we will take the lexicographically smallest $(\mu,n)$-difference set in a set-isotopy class to be canonical.

\section{Algorithms} 

There has been a history of errors in the enumeration of Latin squares (this history is described in \cite{smallLS}). 
As such, it has become standard practice in the enumeration of Latin squares and related structures to enumerate using multiple different methods and check the results are identical.   
We present an algorithm that completed the enumeration of $(\mu,n)$-difference sets that were unique up to set-isotopism for $\mu \geq 2$ and $n \leq 18$, which includes information about $\text{Is}_l(\mathcal{L})$, $\text{Red}_s(\mathcal{L})$, and $\text{Red}_l(\mathcal{L})$ for each $(\mu,n)$-difference set $\mathcal{L}$. 
This information can be used to deduce the number of list-reduced $(\mu,n)$-difference sets. The program also includes a second method to count the number of list-reduced $(\mu,n)$-difference sets by another method, meaning we have a way of verifying the accuracy of our results. 

Within our program we represent a $(\mu,n)$-difference set by the list of columns $(C_1, \ldots, C_{\mu})$, where $(C_{1} (i), \ldots, C_{\mu} (i))$ is an element of our difference set for each $1 \leq i \leq n$. 
Our algorithm will augment a further index to each element of the $(\mu,n)$-difference set by augmenting a column $C$ to the stored list of columns $(C_1, \ldots, C_{\mu})$. 
Such a column, say $C_{\mu+1}$, has to satisfy the necessary and sufficient conditions:
1/ each symbol in $\{0, \ldots, n-1\}$ occurs in $C_{\mu+1}$; 
2/ the column $C_{\mu+1}$ is disjoint to all other columns (i.e. $C_{\mu+1}(k)\neq C_i(k)$); and  
3/ the column $C_{\mu+1}$ is nearly orthogonal with all other columns (i.e. $\{C_{\mu+1}(k)-C_i(k)\} =\{1, \ldots, n/2,n/2, \ldots, n-1\}$). 

Let $R=\{r_1, \ldots, r_n\}$, $S=\{\sigma_1, \ldots, \sigma_n\}$, and $D_i=\{d_{i0}, \ldots, d_{i(n-1)}\}$ be pairwise disjoint sets, representing the rows, symbols, and differences of a column we would want to add, with $1 \leq i \leq \mu$. 
For each cell of $C$ in row $r_j$ containing symbol $\sigma_k$, consider the set $\{r_j, \sigma_k, d_{1 (C_{1}(j) - k)}\}$. 
The collection of such sets, $\mathcal{F}$, when we consider all cells in some potential column must be a partition of the multiset $X = R \cup S \cup \bigcup\limits_{i=1}^{\mu} \{d_{i1}, \ldots, d_{i\frac{n}{2}},d_{i\frac{n}{2}}, \ldots, d_{i(n-1)}\}$. 
As such, this is an instance of the exact cover problem, and we can therefore use existing software such as \textsc{libexact} to find all possible columns that may be augmented to our list $(C_1, \ldots, C_{\mu})$ to give a $(\mu+1,n)$-difference set.

We may assume the first column is in natural order, and refer to such a column as $I$. 
The algorithm begins by generating all columns that could be augmented to this first column using \textsc{libexact}. 
If these $(2,n)$-difference sets are canonical, we place the second column into \textsc{list1}. 
For each second column $C$ in \textsc{list1} we again use \textsc{libexact} to find all columns that may successfully be augmented to $(I,C)$, and we place all such third columns into \textsc{list2}. 
Construct a graph with vertices in \textsc{list2}, and edges connecting points $C_a$ and $C_b$ if $(C_a,C_b)$ is a $(\mu,n)$-difference set. 
Then each clique $(e_1, \ldots, e_{\alpha})$  yields a $(\alpha+2,n)$-difference set, given by $(I,C,e_1, \ldots, e_{\alpha})$. 
For each clique, if the generated $(\alpha+2,n)$-difference set, $\mathcal{L}$, is set-canonical, we calculate $ |\text{Red}_s(\mathcal{L})|$, $ |\text{Is}_l(\mathcal{L})|$, and $|\text{Is}_s(\mathcal{L})|$ and store this information. 
After completion, we merge the results and use Lemma \ref{numberClassesPerClass}  to find the total number of classes. 

Finding cliques is usually a hard problem. 
This is not an issue for our calculations as the clique size of our problem turns out to be very small. 
In fact, no cliques of size three existed in our graph, and the computation time to prove this was  negligible within our program as a whole. 

We used a brute force search over the set-isotopism operators to test canonicality and to calculate $ |\text{Red}_s(\mathcal{L})|$, $ |\text{Is}_l(\mathcal{L})|$, and $|\text{Is}_s(\mathcal{L})|$. Again, the time required for this operation was negligible. 

We also performed our checking step at the point of finding the clique. Each $C_1$ in \textsc{list1} corresponds to $2 \phi(n) n /|\text{Is}_s((I,C_1))|$ list-reduced $(2,n)$-difference sets. Each clique of size $\mu-2$ has $(\mu-2)!$ ways of being augmented onto the list $(I,C_1)$, yielding $(\mu-2)! 2 \phi(n) n /|\text{Is}_s((I,C_1))|$ list-reduced $(\mu,n)$-difference sets. Summing over each of these values gives the total number of list-reduced $(\mu,n)$-difference sets.

\begin{algorithm}[htp]
    \SetKwInOut{Input}{input}
    \SetKwInOut{Output}{output}
\Input{ 
$n$;}

\Output{A list with entries being lists of four integers $(\mathrm{Is}_s,\mathrm{Is}_l, \mathrm{Red}_s, count)$, where count is the number of set-isotopy classes of $(\mu,n)$-difference sets $\mathcal{L}$ that have 
$(|\mathrm{Is}_s(\mathcal{L})|,|\mathrm{Is}_l(\mathcal{L})|, |\mathrm{Red}_s(\mathcal{L})|) = (\mathrm{Is}_s,\mathrm{Is}_l, \mathrm{Red}_s)$
\;}
 
 $store \gets \emptyset$\;
 $\textsc{list1} \gets \emptyset$\;

 \For{$C_1 \in exactcover(I)$}{
 $\textsc{list1} \gets \textsc{list1} \cup C_1$\\
 }
 \For{$C_1 \in \textsc{list1}$}{
 
 $\mathrm{vert}(graph) \gets \emptyset$\;
 $\mathrm{edge}(graph) \gets \emptyset$\;
 
 $\textsc{list2} \gets \emptyset$\;
  \For{$C_2 \in exactcover(I,C_1)$}{
 $\textsc{list2} \gets \textsc{list2} \cup C_2$\\
 $\mathrm{vert}(graph) \gets \mathrm{vert}(graph)\cup C_2$\\
 }

 \For{$v_1 \in \textsc{list2}$}{
 
 \For{$v_2 \in \textsc{list2}$}{
 \If{$(v_1,v_2)$ are a $(\mu,n)$-difference set}{
$edge \gets edge \cup \{\{v_1,v_2\}\}$\\
}
 }
 }
\For{all cliques $(\alpha_1, \ldots, \alpha_{\mu-2})$ of size $\mu-2$ such that $(I, C_1, \alpha_1, \ldots, \alpha_{{\mu}-2})$ is set-canonical}{
	{
	$\mathcal{L} \gets$ the $(\mu,n)$-difference set $(I, C_1, \alpha_1, \ldots, \alpha_{\mu-2})$ \\
	{$store \gets store.add(\mathrm{Is}_s(\mathcal{L}), \mathrm{Is}_l(\mathcal{L}), \mathrm{Red}_s(\mathcal{L}))$
}	
	   }
}

}

 \caption{Algorithm}
\end{algorithm}

\section{Results and conclusions}

The counts that were found appear in Tables \ref{tab:setIso}, \ref{tab:setRed}, \ref{tab:listIso}, and \ref{tab:listRed}. 
Comparing these results to the previously known cases in Table \ref{tab:NOLS}, we see that the new values of particular significance  are when $\mu=3$ and $n \in \{10,12,14,16,18\}$, when $\mu  = 4$ and $n \in \{12,14,16,18\}$, and when $\mu = 5$ and $n \in \{ 14,16,18\}$. The results when $\mu =5$ disproves Conjecture 5.2 of  \cite{NOLS} that proposed the maximum $\mu$ for which a set of $\mu$ cyclic  $\text{MNOLS}$ of order $n$ exists is $\lceil n/4\rceil +1$, as there does not exist a set of five $\text{MNOLS}$ of order $14$, $16$, and $18$ as predicted by the conjecture.(In fact, the conjecture predicted a set of $6$ cyclic $\text{MNOLS}$ of order $18$.)

For $n =14$ the search consumed $14$ minutes of CPU time and for $n =16$ the search consumed $77$ hours of CPU time, with negligible amounts of memory. 
As a means of comparison, running a depth first search for set-isotopic classes consumed $20$ hours and $154$ days of CPU time, respectively, and required $7$GB of RAM (to save on repeated computation). The resulting $(\mu,n)$-difference sets are of reasonable size, and are provided online for examination and verification \cite{Store_Dropbox}.

For $n=18$, the search was conducted independently by two authors, consuming $2337$ core-days and $4137$ core-days, each with memory usage of only a few MB\footnote{The times reported here refer to a logical core using hyperthreading; with a single
physical core per thread the core-time would be around $25\%$ smaller}. 
With this program, generating the sets of $2\text{MNOLS}$ of order $20$ would take
approximately $50$ core-hours, and there are around $1.1\cdot 10^9$ set-isotopy classes.
For each of those, finding all $\mu$-MNOLS would take approximately $50$ core-seconds on average, so the total running time for $n = 20$ would be around
$2000$ core-years.

\begin{table}[!h]
\begin{center}
\resizebox{\textwidth}{!}{%
  \begin{tabular}{ |c | c | c|c | c | c|c | c | c|c | c | c|c | c | c|c | c | c|c|  }
    \hline
$n$& $4$& $6$& $8$& $10$& $12$& $14$&$16$&$18$ \\ \hhline{|=|=|=|=|=|=|=|=|=|}  

$\mu=2$&1&2&9&68&1140&19040&489296&28303688\\ \hline 			
$\mu=3$&0&1&1&73&4398& 429111&70608753&31992833620\\  \hline 			
$\mu=4$&0&0&0&1&2&117&14672&8354783\\ \hline 		
$\mu=5$&0&0&0&0&0&0&0&0\\ \hline 		
\end{tabular}
}
\caption{The number of collections of $\mu$ cyclic $\text{MNOLS}$ of order $n$ under set-isotopy equivalence.}\label{tab:setIso}
\end{center}
\end{table}

\begin{table}[!h]
\begin{center}
\resizebox{\textwidth}{!}{%
  \begin{tabular}{ |c | c | c|c | c | c|c | c | c|c | c | c|c | c | c|c | c | c|c|  }
    \hline
$n$& $4$& $6$& $8$& $10$& $12$& $14$&$16$&$18$ \\ \hhline{|=|=|=|=|=|=|=|=|=|}  

$\mu=2$&2&12&136&2340&52608&1589056&62516224&3056224608\\ \hline
$\mu=3$&0&6&16&2920&211104& 36031716&9037728896&3455226014904\\ \hline 		
$\mu=4$&0&0&0&20&96&8638&1870592&902182968\\ \hline  	
$\mu=5$&0&0&0&0&0&0&0&0\\ \hline 		

\end{tabular}
}
\caption{The number of collections of $\mu$ cyclic $\text{MNOLS}$ of order $n$ under set-reduced equivalence.}\label{tab:setRed}
\end{center}
\end{table}

\begin{table}[!h]
\begin{center}
\resizebox{\textwidth}{!}{%
  \begin{tabular}{ |c | c | c|c | c | c|c | c | c|c | c | c|c | c | c|c | c | c|c|  }
    \hline
$n$& $4$& $6$& $8$& $10$& $12$& $14$&$16$&$18$ \\ \hhline{|=|=|=|=|=|=|=|=|=|}  

$\mu=2$&1&3&12&128&2224&38000&977696&56603408\\ \hline
$\mu=3$&0&2&6&438&26388& 2574306&423652518&191957000556\\ \hline 			
$\mu=4$&0&0&0&12&48&2484&350730&200481924\\ \hline 	
$\mu=5$&0&0&0&0&0&0&0&0\\ \hline 		

\end{tabular}
}
\caption{The number of collections of $\mu$ cyclic $\text{MNOLS}$ of order $n$ under list-isotopy equivalence.}\label{tab:listIso}
\end{center}
\end{table}

\begin{table}[!h]
\begin{center}
\resizebox{\textwidth}{!}{%
  \begin{tabular}{ |c | c | c|c | c | c|c | c | c|c | c | c|c | c | c|c | c | c|c|  }
    \hline
$n$& $4$& $6$& $8$& $10$& $12$& $14$&$16$&$18$ \\ \hhline{|=|=|=|=|=|=|=|=|=|}  %

$\mu=2$&4&24&256&4640&105216&3178112& 125026304&6112406016 \\ \hline 		
$\mu=3$&0&12&96&17520&1266624& 216190296&54226373376&20731356060048\\ \hline 		
$\mu=4$ &0&0&0&480&2304&207312&44879616&21652047792 \\ \hline 		
$\mu=5$&0&0&0&0&0&0&0&0\\ \hline 		
\end{tabular}
}
\caption{The number of collections of $\mu$ cyclic $\text{MNOLS}$ of order $n$ under list-reduced equivalence.}\label{tab:listRed}
\end{center}
\end{table}

We say a list of $\mu$ cyclic $\text{MNOLS}$ of order $n$ is of \emph{type 0} if it is isotopically equivalent to a list of reduced  $\mu$ cyclic $\text{MNOLS}$ of order $n$, $\mathcal{L} = (L_1, \ldots, L_\mu)$, with $(0,0,1),(1,0,0) \in L_2$, and is of \emph{type 1} otherwise. 
A set of $\mu$ $\text{MNOLS}$ of order $n$ is of type 0 if fixing the order in some way gives a list of $\mu$ $\text{MNOLS}$ of order $n$ of type 0, and is of type 1 otherwise.

A collection of $\mu$ cyclic  $\text{MNOLS}$ of order $n$ contains a \emph{row-intercalate} of difference $d$ if two of its Latin squares $L$ and $M$ have two symbols $e,e'$ with $e<e'$ and $e'-e=d$ such that $L(r,0) = M(r',0)=e$ and also $L(r',0) = M(r,0) = e'$, for some $r,r' \in \{0,\ldots, n-1\}$. 
Then it is clear that a collection of $\mu$ cyclic $\text{MNOLS}$ of order $n$ is of type 0 if and only if it contains a row-intercalate of difference $d$ and $\gcd(d,n)=1$. 
Clearly set-isotopy preserves type.
In Tables \ref{tab:tab2wayTypes}, \ref{tab:tab3wayTypes}, \ref{tab:tab4wayTypes}, \ref{tab:tab2wayTypes16}, \ref{tab:tab3wayTypes16}, and \ref{tab:tab4wayTypes16} we show the number of set-isotopy classes of each type. 
Observe that the proportion of set-isotopy classes that are of type 0 increases as $\mu$ increases. 
This may be of interest in future searches for sets of $\mu$ cyclic $\text{MNOLS}$ of order $n$ where $\mu$ is relatively large. 
Considering each type individually may allow more efficient construction of those set-isotopy classes with non-trivial set-autotopy group, as each set-autotopy must map row-intercalates to row-intercalates. 
Note that $|\text{Red}_s(\mathcal{L})|=1$ for $n=14$, so we omit the column for $|\text{Red}_s(\mathcal{L})|$ in this case.

\begin{table}[!h]
\begin{center}
  \begin{tabular}{ |c | c | c|c | c | c|c | c | c|c | c | c|c | c | c|c | c | c|c|  }
    \hline
$|\text{Is}_s(\mathcal{L})|$ & $|\text{Is}_l(\mathcal{L})|$ & \#Type 0&\#Type 1&\#Total  \\ \hhline{|=|=|=|=|=|=|=|=|}
1&1&3618&15186&18804\\ \hline 		
2&1&0&80&80\\ \hline
2&2&46&88&134\\ \hline
3&3&2&14&16\\ \hline
6&6&1&5 &6\\ \hline
&total:&3667&15373 &19040\\ \hline
\end{tabular}
\caption{The number of collections of two cyclic $\text{MNOLS}$ of order $14$, by their type and autotopy group sizes.} \label{tab:tab2wayTypes}
\end{center}
\end{table}

\begin{table}[!h]
\begin{center}
  \begin{tabular}{ |c | c | c|c | c | c|c | c | c|c | c | c|c | c | c|c | c | c|c|  }
    \hline
$|\text{Is}_s(\mathcal{L})|$ & $|\text{Is}_l(\mathcal{L})|$ & \#Type 0 & \#Type 1 & \#Total  \\ \hhline{|=|=|=|=|=|=|=|=|}
1&1&202382&226436&428818\\ \hline 				
2&2&146&57&203\\ \hline
3&1&24&63&87\\ \hline
6&2&1&2&3 \\ \hline
&total:&202553&226558&429111\\ \hline

\end{tabular}
\caption{The number of collections of three cyclic $\text{MNOLS}$ of order $14$, by their type and  autotopy group sizes.} \label{tab:tab3wayTypes}
\end{center}
\end{table}

\begin{table}[!h]
\begin{center}
  \begin{tabular}{ |c | c | c|c | c | c|c | c | c|c | c | c|c | c | c|c | c | c|c|  }
    \hline
$|\text{Is}_s(\mathcal{L})|$ & $|\text{Is}_l(\mathcal{L})|$ &\#Type 0 & \#Type 1 & \#Total   \\ \hhline{|=|=|=|=|=|=|=|=|}

1&1&67&26&93\\ \hline 			
2&1&3&8&11\\ \hline
2&2&1&0&1\\ \hline
3&1&4&7&11\\ \hline
6&2&1&0 &1\\ \hline
&total:&76&41&117\\ \hline

\end{tabular}
\caption{The number of collections of four cyclic $\text{MNOLS}$ of order $14$, by their type  and autotopy group sizes.} \label{tab:tab4wayTypes}
\end{center}
\end{table}

\begin{table}[!h]
\begin{center}
  \begin{tabular}{ |c | c | c|c | c | c|c | c | c|c | c | c|c | c | c|c | c | c|c|  }
    \hline
$|\text{Is}_s(\mathcal{L})|$ & $|\text{Is}_l(\mathcal{L})|$ & $|\text{Red}_s(\mathcal{L})|$  & \#Type 0&\#Type 1&\#Total  \\ \hhline{|=|=|=|=|=|=|=|=|}

1& 1& 1& 106794& 380686& 487480 \\ \hline
2& 1& 1& 12& 822& 834\\ \hline
2& 2& 1& 260& 660& 920\\ \hline
4& 2& 1& 0& 12& 12\\ \hline
2& 1& 2& 46& 0& 46\\ \hline
4& 2& 2& 4& 0& 4\\ \hline
&&total:&107116&382180&489296\\ \hline
\end{tabular}
\caption{The number of collections of two cyclic $\text{MNOLS}$ of order $16$, by their type and autotopy group sizes.} \label{tab:tab2wayTypes16}
\end{center}
\end{table}

\begin{table}[!h]
\begin{center}
  \begin{tabular}{ |c | c | c|c | c | c|c | c | c|c | c | c|c | c | c|c | c | c|c|  }
    \hline
$|\text{Is}_s(\mathcal{L})|$ & $|\text{Is}_l(\mathcal{L})|$ & $|\text{Red}_s(\mathcal{L})|$  & \#Type 0&\#Type 1&\#Total  \\ \hhline{|=|=|=|=|=|=|=|=|}
1& 1& 1& 36845488& 33760273& 70605761 \\ \hline
2& 2& 1& 2326& 666& 2992\\ \hline
&&total:&36847814&33760939&70608753\\ \hline
\end{tabular}
\caption{The number of collections of three cyclic $\text{MNOLS}$ of order $16$, by their type and autotopy group sizes.} \label{tab:tab3wayTypes16}
\end{center}
\end{table}

\begin{table}[!h]
\begin{center}
  \begin{tabular}{ |c | c | c|c | c | c|c | c | c|c | c | c|c | c | c|c | c | c|c|  }
    \hline
$|\text{Is}_s(\mathcal{L})|$ & $|\text{Is}_l(\mathcal{L})|$ & $|\text{Red}_s(\mathcal{L})|$  & \#Type 0&\#Type 1&\#Total  \\ \hhline{|=|=|=|=|=|=|=|=|}
1& 1& 1& 11146& 3401& 14547 \\ \hline
2& 1& 1& 28& 79& 107\\ \hline
2& 2& 1& 7& 2& 9\\ \hline
2& 1& 2& 8& 0& 8\\ \hline
4& 1& 4& 1& 0& 1\\ \hline
&&total:&11190&3482&14672\\ \hline
\end{tabular}
\caption{The number of collections of four cyclic  $\text{MNOLS}$ of order $16$, by their type and autotopy group sizes.} \label{tab:tab4wayTypes16}
\end{center}
\end{table}

\begin{table}[!h]
\begin{center}
  \begin{tabular}{ |c | c | c|c | c | c|c | c | c|c | c | c|c | c | c|c | c | c|c|  }
    \hline
$|\text{Is}_s(\mathcal{L})|$ & $|\text{Is}_l(\mathcal{L})|$ & $|\text{Red}_s(\mathcal{L})|$  & \#Type 0&\#Type 1&\#Total  \\ \hhline{|=|=|=|=|=|=|=|=|}

1& 1& 1& 4378529 &23914135 &28292664\\ \hline
2 &1 &1 &0 &3568& 3568\\ \hline
2 &2& 1 &1642 &5414 &7056\\ \hline
2 &1 &2 &400& 0 &400\\ \hline
&&total:&4380571& 23923117&28303688\\ \hline
\end{tabular}
\caption{The number of collections of two cyclic $\text{MNOLS}$ of order $18$, by their type and autotopy group sizes.} \label{tab:tab2wayTypes16}
\end{center}
\end{table}

\begin{table}[!h]
\begin{center}
  \begin{tabular}{ |c | c | c|c | c | c|c | c | c|c | c | c|c | c | c|c | c | c|c|  }
    \hline
$|\text{Is}_s(\mathcal{L})|$ & $|\text{Is}_l(\mathcal{L})|$ & $|\text{Red}_s(\mathcal{L})|$  & \#Type 0&\#Type 1&\#Total  \\ \hhline{|=|=|=|=|=|=|=|=|}
1& 1 &1& 12650027871& 19342805458 &31992833329 \\ \hline
3& 1 &1 &47& 176& 223\\ \hline
3 &1& 3 &0& 68& 68\\ \hline
&&total: &12650027918 &19342805702 &31992833620\\ \hline
\end{tabular}
\caption{The number of collections of three cyclic $\text{MNOLS}$ of order $18$, by their type and autotopy group sizes.} \label{tab:tab3wayTypes16}
\end{center}
\end{table}

\begin{table}[!h]
\begin{center}
  \begin{tabular}{ |c | c | c|c | c | c|c | c | c|c | c | c|c | c | c|c | c | c|c|  }
    \hline
$|\text{Is}_s(\mathcal{L})|$ & $|\text{Is}_l(\mathcal{L})|$ & $|\text{Red}_s(\mathcal{L})|$  & \#Type 0&\#Type 1&\#Total  \\ \hhline{|=|=|=|=|=|=|=|=|}
1& 1 &1& 5291250& 3060794 &8352044 \\ \hline
2& 1 &1 &675 &1799& 2474\\ \hline
2& 1& 2 &265& 0& 265\\ \hline
&&total:& 5292190 &3062593 &8354783\\ \hline
\end{tabular}
\caption{The number of collections of four cyclic  $\text{MNOLS}$ of order $18$, by their type and autotopy group sizes.} \label{tab:tab4wayTypes16}
\end{center}
\end{table}




\ 

\

\section*{References}

\bibliographystyle{jcdStyle1}


\begin{thebibliography}{10}
\expandafter\ifx\csname urlstyle\endcsname\relax
  \providecommand{\doi}[1]{doi:\discretionary{}{}{}#1}\else
  \providecommand{\doi}{doi:\discretionary{}{}{}\begingroup
  \urlstyle{rm}\Url}\fi





\bibitem{PermGroupsCameron}
P.~Cameron.
\newblock \emph{Permutation Groups}.
\newblock Cambridge University Press, Cambridge, England, 1999.

\bibitem{DCA}
F.~Demirkale, D.M. Donovan, J.L. Hall, A.~Khodkar, and A.~Rao, Difference
  covering arrays and pseudo-orthogonal Latin squares, \emph{Graphs and
  Combinatorics} 32 no.4 (2016), 1353--1374.

\bibitem{DirConst}
F.~Demirkale, D.M. Donovan, and A.~Khodkar, Direct constructions for general
  families of cyclic mutually nearly orthogonal latin squares, \emph{J. Combin.
  Des.} 23 no.5 (2015), 195--203.

\bibitem{Store_Dropbox}
F.~Demirkale, D.M. Donovan, J.I.~Kokkala, and T.G.~Marbach
\newblock \emph{MNOLS data}.
\newblock https://tinyurl.com/MNOLSenum.

\bibitem{EnumMOLS}
J.~Egan, and I.M. Wanless, Enumeration of {MOLS} of small order, \emph{Math.
  Comp.} 85 no.298 (2016), 799--824.

\bibitem{Faradzev}
I.A. Farad\v{z}ev, Generation of nonisomorphic graphs with a given distribution
  of the degrees of vertices. (russian), \emph{Algorithmic studies in
  combinatorics (Russian)} 185 (1978), 11--19.

\bibitem{LS11}
A.~Hulpke, P.~Kaski, and P.R.J. \"{O}sterg\r{a}rd, The number of latin squares
  of order 11, \emph{Math. Comp.} 80 no.274 (2011), 1197--1219.

\bibitem{ClassAlgorithms}
P.~Kaski, and P.R.J. \"{O}sterg\r{a}rd.
\newblock \emph{Classification algorithms for codes and designs}.
\newblock Algorithms {C}omput. {M}ath. 15, {S}pringer (2006).

\bibitem{LS8}
G.~Kolesova, C.W.H. Lam, and L.~Thiel, On the number of $8\times 8$ latin
  squares, \emph{J. Combin. Theory Ser. A} 54 no.1 (1990), 143--148.

\bibitem{DisMathUsingLS}
C.~Laywine, and G.~Mullen.
\newblock \emph{Discrete mathematics using latin squares}.
\newblock Wiley-Interscience Series in Discrete Mathematics and Optimization
  (1998).

\bibitem{NOLS}
P.C. Li, and G.H.J van Rees, Nearly orthogonal latin squares, \emph{J. Combin.
  Math. Combin. Comput.} 62 (2007), 13--24.

\bibitem{IsoFreeExhaGen}
B.D. McKay, Isomorph-free exhaustive generation, \emph{J. Algorithms} 26 no.2
  (1998), 306--324.

\bibitem{smallLS}
B.D. McKay, A.~Meynert, and W.J. Myrvold, Small latin squares, quasigroups and
  loops, \emph{J. Combin. Des.} 15 no.2 (2007), 98--119.

\bibitem{nauty}
B.D. McKay, and A.~Piperno, Practical graph isomorphism, {II}, \emph{J.
  Symbolic Comput.} 60 (2014), 94--112.

\bibitem{MostLSsubsquares}
B.D. McKay, and I.M. Wanless, Most latin squares have many subsquares, \emph{J.
  Combin. Theory Ser. A} 86 no.2 (1999), 322--347.

\bibitem{numberLS}
B.D. McKay, and I.M. Wanless, On the number of latin squares, \emph{Ann. Comb.}
  9 no.3 (2005), 335--344.

\bibitem{no4MNOLS6}
E.B. Pasles, and D.~Raghavarao, Mutually nearly orthogonal latin squares of
  order 6, \emph{Util. Math.} 65 (2004), 65--72.

\bibitem{MNOLSinvent}
D.~Raghavarao, S.S. Shrikhande, and M.S. Shrikhande, Incidence matrices and
  inequalities for combinatorial designs, \emph{J. Combin. Des.} 10 no.1
  (2002), 17--26.

\bibitem{EveryOne}
R.C. Read, Every one a winner or how to avoid isomorphism search when
  cataloguing combinatorial configurations., \emph{Algorithmic aspects of
  combinatorics (Conf., Vancouver Island, B.C., 1976). Ann. Discrete Math.} 2,
  1978, pp.\ 107--120.

\bibitem{PrivateVanRees}
G.H.J. van Rees.
\newblock Private communication (2015).

\end{thebibliography}
\end{document}